\documentclass[Afour,sageh,times]{sagej}

\usepackage[colorlinks,bookmarksopen,bookmarksnumbered,citecolor=red,hypertexnames=false]{hyperref}
\usepackage{breakurl}

\def\volumeyear{2019}

\usepackage{color}

\usepackage[T1]{fontenc}

\usepackage{graphics}

\usepackage{amsmath,amssymb}

\usepackage[]{algorithm2e}
\usepackage{verbatim}

\usepackage{subcaption}

\usepackage[english]{babel}

\newcommand{\fv}{\mathbf f}
\newcommand{\ev}{\mathbf e}

\newcommand{\uv}{\mathbf{u}}

\newcommand{\nv}{\mathbf n}

\newcommand{\Ev}{ \mathbf{E}}
\newcommand{\nullv}{\mathbf{0}}
\newcommand{\Pv}{ \mathbf{P}}

\newcommand{\Dv}{ \mathbf{D}}

\newcommand{\sigmat}{\boldsymbol{\sigma}}

\newcommand{\epst}{\boldsymbol{\varepsilon}}
\newcommand{\epsilont}{\boldsymbol{\epsilon}}
\newcommand{\betat}{\boldsymbol{\beta}}

\newcommand{\St}{\mathbf{S}}
\newcommand{\Ct}{\mathbf{C}}
\newcommand{\dt}{\mathbf{d}}

\newcommand{\It}{\mathbf{I}}
\newcommand{\Et}{\mathbf{E}}

\newcommand{\nullt}{\mathbf{0}}

\newcommand{\opdiv}{\operatorname{div}}
\newcommand{\opd}{ {\operatorname{d}} }

\newcommand{\rev}[1]{ { #1} }

\newcommand{\SI}[2]{\ensuremath{#1\hspace{1pt}\text{#2}}}

\def\meter{m}
\def\milli{m}

\begin{document}
	
	\graphicspath{
		{./}{./pics/}{./picsPs/} }

\runninghead{Meindlhumer, Pechstein and Humer}

\title{Variational inequalities for ferroelectric constitutive modeling}

\author{Martin Meindlhumer \affilnum{1}, Astrid Pechstein\affilnum{1} and Alexander Humer\affilnum{1}}

\affiliation{\affilnum{1}Johannes Kepler University Linz, Austria}

\corrauth{Institute of Technical Mechanics,
	Johannes Kepler Universtity Linz,
	Altenberger Str. 69,
	4040 Linz,
	AUSTRIA.}

\email{martin.meindlhumer@jku.at}

\begin{abstract}
This paper is concerned with modeling the polarization process in ferroelectric media.
We develop a thermodynamically consistent model, based on phenomenological descriptions of free energy as well as switching and saturation conditions in form of inequalities. 
Thermodynamically consistent models naturally lead to variational formulations. We propose to use the concept of variational inequalities.
We aim at combining the different phenomenological conditions into one variational inequality.
In our formulation we use
one Lagrange multiplier for each condition (the onset of domain switching and saturation), each
satisfying Karush-Kuhn-Tucker conditions.
An update for reversible and remanent quantities is then computed within one, in general nonlinear, 
iteration. 
\end{abstract}


\maketitle
\section{Introduction}

In the current work we aim at describing the process of polarziation in ferroelectric media.
Polarization has been described as a dissipative process in a thermodynamically consistent
framework based on the Helmholtz free energy in a series of papers by 
\cite{BassiounyGhalebMaugin:1988a,BassiounyGhalebMaugin:1988b,BassiounyMaugin:1989a,BassiounyMaugin:1989b}.
The notions introduced there are close to the theory of elasto-plasticity, including
internal variables, yield conditions and hardening moduli. The evolution of the internal
variables is based on inequality constraints, such as the so-called switching condition.
This condition defines the onset of remanent polarization in ferroelectric media.

\cite{McMeekingLandis:2002} as well as \cite{Landis:2002} 
developed their theory of multi-axial polarization based on those early works.
The internal variables are the polarization vector in \cite{McMeekingLandis:2002}, and
an independent polarization strain is added in \cite{Landis:2002}.
They specified non-linear free energy functions that account also for the saturation
phenomenon. Once reaching saturation, the polarization cannot grow any further, but
may change direction or be reduced by (electric or mechanic) depolarization.

\cite{KamlahTsakmakis:1999} chose a different set of conditions to
describe the polarization process, see also the overview given in \cite{Kamlah:2001}. They
assume the reversible part of the free energy to be quadratic, and add separate switching
and saturation conditions for remanent polarization and remanent polarization strain.

All of these frameworks have been implemented in finite element codes by various authors.
\cite{KamlahBoehle:2001} provided the increments for the evolution polarization and
polarization strain for the material described in \cite{KamlahTsakmakis:1999}.
In \cite{Klinkel:2006}  a return mapping algorithm is presented realizing Landis-type
material properties. \cite{SemenovLiskowskyBalke:2010} provide such an
algorithm for a finite element formulation including the vector (not scalar) potential
of the dielectric displacement vector. 
\cite{ElhadrouzZinebPatoor:2005} implemented a material model including different switching functions for polarization
and polarization strain, similar as proposed by \cite{KamlahTsakmakis:1999}. 
\cite{ZouariZinebBenjeddou:2011} ~proposed a quadratic element realizing a constitutive law, similar to the one proposed by \cite{Kamlah:2001}.

\rev{A hybrid phenomenological model for ferroelectric ceramics is presented by \cite{stark2016hybrid,stark2016hybrid1}.  Comments on the stability of ferroelectric constitutive models were given by  \cite{StarkNeumeisterBalke:2016a} and \cite{BotteroIdiart:2018}. A general mathematical framework allowing for proofs of existence and uniqueness of solutions  was recently presented  by \cite{Pechstein2019Polarization}.  \cite{humer2020nonlinear} presented a formulation considering large deformation and hysteresis effects using threefold multiplicative decomposition of the deformation gradient tensor into a reversible mechanical, a reversible electrical and an irreversible part. }  

We will follow Kamlah's approach, but formulate the theory in the framework of variational
inequalities. These inequalities arise naturally from energy-based constitutive models
with dissipation, see \cite{MieheRosatoKiefer:2011} for the application to ferroelectric
and also ferromagnetic materials. In our case, we add both switching and saturation
condition. In the following derivations, we will see that,
in the variational framework, these conditions are not independent as proposed in \cite{Kamlah:2001},
but that the Lagrangian multiplier enforcing saturation becomes an additional term in the switching
criterion. This way, both criteria can be satisfied at once. 

\rev{
	The framework of variational inequalities allows the direct implementation of different inequality constraints in a thermodynamically consistent way. In contrast to return mapping algorithms the equations may be solved within one in general nonlinear iteration. 
	An alternative approach to model the evaluation of internal variables is the construction of nonlinear  dissipation functions as proposed by \cite{Landis:2002}. An advantage of the latter method is, that it can be easily implemented as no inequalities have to be taken into account. A drawback of these dissipation functions is, that they are in general not differentiable.   }

Concerning the finite element implementation, we use mixed finite elements for the mechanical
unknowns, adding independent stresses. We use non-standard TDNNS elements, where tangential displacements
and the normal component of the normal stress are chosen as degrees of freedom \cite{PechsteinSchoeberl:2011,PechsteinSchoeberl:2016}.
These elements were adapted to linear piezoelectricity in \cite{PechsteinMeindlhumerHumer:2018,MeindlhumerPechstein:2018}. 
\rev{The TDNNS elements have been shown to be free from locking effects and therefore they are highly suitable for the efficient discretization of thin structures, which will be shown in one particular example.  }

Our contribution is organized as following: In the first section we provide a phenomenological material law in the spirit of \cite{Kamlah:2001}. In the next section we show how this material law can be embedded into incremental variational formulations similar to \cite{MieheRosatoKiefer:2011}. As switching and saturation are modeled via inequality constraints, a variational inequality is obtained. We show how these variational inequalities can be implemented in a finite element formulation. In the last section of this contribution we show several numerical results. We consider benchmark problems \rev{from} the literature, including electrical polarization, mechanical depolarization, non-proportional loading and  depolarization under bending stress. Additionally we show the polarization of a bimorph structure. Last, a short review is provided. 

\section{A phenomenological material law} \label{sec:phenomenologic}

In this section we briefly introduce the phenomenological material model in the spirit of  \citet{Kamlah:2001}, before \rev{deriving} an energy-based approach in the next section.
We consider a ferroelectric body in the domain $\Omega \rev{\subset} \mathbb{R}^3$ with boundary $\partial\Omega$. We are concerned with \rev{isothermal}, rate independent and quasi-static deformation and polarization processes. We are interested in finding the displacement field $\uv$, the electric potential $\phi$ and the remanent polarization $\Pv_I$. 

 The electric field is introduced as the negative gradient of the electric potential $\Et=-\nabla\phi$, whereas we use the linear strain tensor $\epst=\frac{1}{2}\left(\nabla\uv+\nabla\uv^T\right)$. 
The mechanical balance equation and Gauss' law are
\begin{eqnarray}
	-\opdiv\sigmat&=&\fv, \label{eq:balance_mech}\\
	-\opdiv\Dv &=& 0, \label{eq:balance_el}
\end{eqnarray}
for the stress tensor $\sigmat$ and the dielectric displacement $\Dv$. The mechanical boundary conditions are
\begin{equation}
\uv = \nullv \text{ on } \Gamma_1 \qquad \text{and} \qquad \sigmat \cdot \mathbf{n} = \sigmat_0 \text { on } \Gamma_2 = \partial \Omega \backslash \Gamma_1,
\end{equation}
and the electrical boundary conditions
\begin{equation}
\phi = \phi_0 \text{ on } \Gamma_3 \qquad \text{and} \qquad \Dv \cdot \mathbf{n} = 0 \text { on } \Gamma_4 = \partial \Omega \backslash \Gamma_3,
\end{equation}
 with $\nv$ denoting the normal vector on the corresponding boundary. The constitutive model relates stress $\sigmat$ and dielectric displacement $\Dv$ to strain $\epst$ and electric field $\Ev$ as well as to the remanent polarization $\Pv_I$. A basic assumption is the additive decomposition of dielectric displacement and elastic strain into a reversible and a remanent or irreversible part. For the reversible parts, constitutive equations analogous to Voigt's linear theory are assumed
 \begin{eqnarray}\label{eq:const1}
 \epst_R=&\epst -\epst_I &= \St^E\cdot\sigmat+\dt^T\cdot\Ev\rev{,}\\
 \Dv_R=&\Dv-\Pv_I &= \dt\cdot\sigmat+\epsilont^\sigma\cdot\Ev.	\label{eq:const2}
 \end{eqnarray}
  The mechanical compliance tensor $\St^E$ as well as the dielectric tensor $ \epsilont^\sigma$ are assumed to be isotropic and constant.
 The components of the \rev{piezoelectric} tensor \rev{$\dt$} depend on $\Pv_I$ via
 \begin{align}\label{eq:dtens}
 \begin{split} 
 d_{kji}=&\frac{\rev{|\Pv_I|}}{P_{sat}}  \big[ d_pe^P_ie^P_je^P_k   \\
 +&d_n \left(\delta_{ij}-e^P_ie^P_j\right)e^P_k  \\
  +&d_t\frac{1}{2} \left(\left(\delta_{ki}-e^P_ke^P_i\right)e^P_j+\left(\delta_{kj}-e^P_ke^P_j\right)e^P_i\right)\big] \end{split}
 \end{align} %
 with unit vector in direction of polarization 
 \begin{equation}\label{eq:ep}
 \mathbf{e}^P=\frac{\Pv_I}{|\Pv_I|},
 \end{equation}
 and the Kronecker delta 
 \begin{equation}\label{eq:deltaij}
 \delta_{ij} =
 \begin{cases}
 1, &         \text{if } i=j,\\
 0, &         \text{if } i\neq j.
 \end{cases}
 \end{equation}

 The unpolarized material does not show any piezoelectric coupling effects. The dielectric tensor grows (linearly) with the norm of $\Pv_I$. The constants $d_p$, $d_n$ and $d_t$ correspond to the standard parameters $d_{33}$, $d_{31}$ and $d_{15}$, respectively, for the fully poled \rev{state} with polarization in direction of $\ev_3$. As proposed by \cite{McMeekingLandis:2002}, we assume the remanent strain  to be volume-preserving and depend only on the polarization via
 \begin{equation}\label{eq:pol_strain}
 \epst_I=\frac{3}{2}\frac{S_{sat}}{P_{sat}^2}\left(\Pv_I\ \otimes\Pv_I\ -\frac{1}{3} \It ~ |\Pv_I|^2\right),
 \end{equation}
 where $S_{sat}$ characterizes the maximum possible remanent strain. The polarization strain $\epst_I$ is a deviatoric uni-axial strain and does not cause any further remanent strain.
The evolution of the remanent polarization $\Pv_I$ is determined by switching and saturation conditions. Once the electric field is increased above the coercive field strength $E_C$, the material will be polarized irreversibly. The switching condition is a condition on the electric driving force $\hat{\Ev}$ which will be specified in the sequel of this contribution. In the simplest case, the absolute value of $\hat{\Ev}$ may not exceed $E_C$,
	\begin{equation}\label{eq:switch_el}
	f_P(\hat{\Ev})=|\hat{\Ev}|-E_C\leq 0.
	\end{equation}
	
On the other hand, the saturation condition states that the remanent polarization is limited and must not exceed the saturation polarization $P_{sat}$, which resembles the fully polarized state. The condition reads
 \begin{equation}\label{eq:sat}
f_S\left(\Pv_I\right)=|\Pv_I|-P_{sat} \leq 0.
\end{equation}

\section{Incremental optimization formulations }

In this section we show how to embed the phenomenological material law from the previous section into  an incremental optimization formulation for (coupled) irreversible processes in the spirit of \cite{MieheRosatoKiefer:2011}. First\rev{,} we briefly {summarize} the incremental formulation for purely reversible processes. This formulation will then be  extended to nonlinear irreversible processes.

\par 

\subsection{Reversible processes}

The stored energy (or free energy) $\Psi^R$ in the domain $\Omega$ is given by
\begin{equation}\label{eq:psi_int}
	\Psi^R=\int_{\Omega}\psi^R \opd \Omega\rev{,}
\end{equation}
with the (volumetric) energy density function $\psi^R$. For coupled piezoelectric problems the energy density function $\psi^R$ reads
\begin{equation}\label{eq:psiR}
	\psi^R=\frac{1}{2} \epst:\sigmat\left(\epst_R,\Dv_R\right)+\frac{1}{2}\Dv\cdot\Ev(\epst_R,\Dv_R).
\end{equation} 

Now, let us consider a finite time interval $\tau:=t-t_n$, $\epst_n$ and $\Dv_n $ denote strain and dielectric displacement at time $t_n$, respectively. Further we assume body forces, boundary forces and electric potential to be constant in the considered time interval.
The work of external loads $W_{ext}^{\tau}$, splits up into two parts
\begin{equation}
	W_{ext}^{\tau}=\rev{W_{ext}^{\tau, V}+W_{ext}^{\tau,B}},
\end{equation}
first \rev{$W_{ext}^{\tau,V}$}, the external work according to body forces,
\begin{equation}
	\rev{W_{ext}^{\tau,V}}=\int_{t_n}^{t}\int_{\Omega}\fv\cdot\dot{\uv}~\operatorname{d\Omega} ~\operatorname{dt},
\end{equation}
 second \rev{$W_{ext}^{\tau,B}$}, external work according to boundary forces and applied electric potential on boundaries
\begin{equation}
	\begin{split}
		\rev{W_{ext}^{\tau,B}}=&\int_{t_n}^{t}\int_{\Gamma_2} \sigmat_0\cdot\dot{\uv}\operatorname{d\Gamma}~\operatorname{dt}\\
		&+\int_{t_n}^{t}\int_{\Gamma_4}\phi_0\dot{\Dv}\cdot\nv\operatorname{d\Gamma}~\operatorname{dt}.
	\end{split}
\end{equation}

  We introduce, according to \cite{MieheRosatoKiefer:2011}, the potential (for reversible processes) $\Pi^{\tau}_R$ with its algorithmic representation
\begin{equation}\label{eq:Pi_tau_R}
	\Pi^{\tau}_R\left(\epst,\Dv\right):=\Psi^R\left(\epst,\Dv\right)-\Psi^R\left(\epst_n,\Dv_n\right) \rev{-} W_{ext}^{\tau}.
\end{equation}
The potential $\Pi^{\tau}_R$ contains the difference between stored energy and work done by external loads $W_{ext}^{\tau}$ in the considered time interval. 
{As the path independence of work is an essential property of linear material behavior, $W_{ext}^{\tau}$ can be determined as the difference of external works between the states at time $t_n$ and $t$.} For given external loads the associated strain and dielectric displacement is given by the constitutive minimization principle 	
\begin{equation}\label{eq:min_Pi_rev}
\begin{split}
	\left(\epst,\Dv\right)=&\operatorname{Arg} \left\{\min_{\epst}\min_{\Dv} \Pi^{\tau}_R\right\}\\
	&\text{s.t.}~\begin{cases}
	\opdiv\Dv=0\rev{,}\\
	\epst=\epst\left(\uv\right).
	\end{cases}
\end{split}		
\end{equation}
The solution of (\ref{eq:min_Pi_rev}) is unique if the potential $\Pi^{\tau}_R$ is convex, which is the case for {stable} piezoelectric materials. 
In the sequel of this contribution we will use, instead of strain $\epst$ and dielectric displacement $\Dv$, mechanical stress $\sigmat$ and electric field $\Ev$ as free variables. The  thermodynamic potential according to this specific choice of free variables is the free enthalpy $H^R$ given in analogy to \eqref{eq:psi_int}, via the enthalpy density $h^R$ by

\begin{equation}\label{eq:H_int}
	H^R=\int_{\Omega}h^R\operatorname{d\Omega}.
\end{equation}
The enthalpy density \rev{$h^R$} is related to the free energy via Legendre transformation, 

\begin{equation}\label{eq:legendre_trans}
	h^R\left(\sigmat,\Ev\right)=\min_{\epst,\Dv}\left(\psi^R\left(\epst,\Dv\right)-\sigmat:\epst-\Dv\cdot\Ev\right).
\end{equation}
For the case of linear piezoelasticity the enthalpy density reads
\begin{equation}
	h^R\left(\sigmat,\Ev\right)=-\frac{1}{2}\sigmat:\St^E:\sigmat-\sigmat:\dt\cdot\Ev-\frac{1}{2}\Ev\cdot\epsilont^{\sigmat}\cdot\Ev.
\end{equation}
In analogy to \eqref{eq:Pi_tau_R} we introduce for a (finite) time interval $\tau=t-t_n$  the potential $\rev{\Pi^{\tau}_H}$ with its algorithmic expression
\begin{equation}\label{eq:Pi_tau_H}
	\rev{\Pi^{\tau}_H}=H^R\left(\sigmat,\Ev\right)-H^R\left(\sigmat_n,\Ev_n\right),
\end{equation}
where suffix $n$ denotes quantities at time $t_n$. Analogously to \eqref{eq:min_Pi_rev} the associated stress and electric field is given by the constitutive maximization principle 
\begin{equation}\label{eq:min_Pi_tau_H}
\begin{split}
	\left(\sigmat,\Ev\right)=&\operatorname{Arg}\left\{\max_{\Ev}\max_{\sigmat}~\rev{\Pi^{\tau}_H}\right\} \\
	&\text{s.t.} 	
	\begin{cases} 
		\opdiv\sigmat=-\fv\rev{,}\\
		\Ev=-\nabla\phi.
	\end{cases}
\end{split}
\end{equation}

Note, that the potential $\rev{\Pi^{\tau}_H}$ in \eqref{eq:Pi_tau_H}, contrary to \rev{$\Pi^{\tau}_R$} in \eqref{eq:Pi_tau_R}  does not contain the external work. External forces and applied potential are rather enforced via the side constraints \rev{and boundary conditions} in \eqref{eq:min_Pi_tau_H}. The first side constraint is the mechanical balance equation, the second ensures Faraday's \rev{law} by introducing the electric field as the negative gradient of the electric potential $\phi$. \par 
The mechanical balance equation is typically treated by introducing the (mechanical) displacement $\uv$ as corresponding Lagrangian multiplier. The enthalpy $H^R$ is extended and reads

\begin{equation}\label{eq:H_ext}
	\mathcal{L}=\int_\Omega h^R-\left(\opdiv\sigmat+\fv\right)\cdot\uv \operatorname{d\Omega}.
\end{equation}
The potential $\rev{\Pi^{\tau}_{\mathcal{L}}}$ corresponding to the extended enthalpy $\mathcal{L}$ is given by
\begin{equation}\label{eq:Pi_tau_L}
	\rev{\Pi^{\tau}_{\mathcal{L}}}=\mathcal{L}\left(\sigmat,\Ev,\uv\right) - \mathcal{L}\left(\sigmat_n,\Ev_n,\uv_n\right).
\end{equation}
The constitutive maximization problem \eqref{eq:min_Pi_tau_H} transforms into a saddle point problem. It reads
\begin{equation}\label{eq:min_}
	\left(\sigmat,\Ev,\uv\right) = \operatorname{Arg}\left\{\max_{\Ev=-\nabla\phi}\max_{\sigmat}\min_{\uv}~\rev{\Pi^{\tau}_{\mathcal{L}}} \right\}.
\end{equation}

\subsection{Dissipative processes}
The material response in ferroelastic materials - for significantly high external loads -  is characterized by non-reversible local response. This contribution focuses on  polarization effects. The process of polarization is a dissipative process. Conducted work is not stored completely as free energy, but some of it is dissipated into heat. 

We introduce the dissipated work $D^{\tau}$ for  the time interval $\tau=t_n-t$ and the considered region $\Omega$. It is given by
\begin{equation}\label{eq:dissipaton}
	D^{\tau}:= \int_{t_n}^{t}\int_\Omega\mathcal{D}\operatorname{d\Omega}~\operatorname{dt}\geq 0,
\end{equation}
with the local dissipation $\mathcal{D}$ .
The second law of thermodynamics states that, for each time interval, $D^{\tau}\geq0$ has to be non-negative. This implies positive dissipation $\mathcal{D}\geq0$ for arbitrary processes. 
In this contribution we follow the approach of internal variables, and  introduce irreversible polarization $\Pv_{I}$. Consequently the stored energy density function for dissipative processes $\psi$ depends also on $\Pv_{I}$,  
\begin{equation}
{\psi}={\psi}\left(\epst,\Dv,\Pv_{I}\right).
\end{equation}
Analogously to \eqref{eq:psi_int} the stored energy in the domain $\Omega$ is given by
\begin{equation}\label{eq:int_Psi_D}
	\Psi=\int_\Omega\psi\operatorname{d\Omega}.
\end{equation}
Note, that in this work, remanent straining is always assumed to depend directly on the remanent polarization as proposed by e.g. \cite{McMeekingLandis:2002} or \cite{LinnemannKlinkel:2009}. Therefore, the irreversible part of the strain $\epst_I$ is not to be considered as a free variable.\par

According to \cite{MieheRosatoKiefer:2011} we introduce the driving force $\hat{\Ev}$. It is defined as the negative derivative of the stored energy density $\psi$ by the internal variables, in our case the irreversible polarization $\Pv_{I}$ 
\begin{equation}\label{eq:driving_force}
	\hat{\Ev}:=-\frac{\partial}{\partial{\Pv_{I}}}\rev{\psi}.
\end{equation}
The driving force $\hat{\Ev}$ is also referred to \rev{as} internal constitutive force, as  $\hat{\Ev}$ is the work conjugate of $\Pv_{I}$. As shown in \cite{MieheRosatoKiefer:2011}, the dissipation is related to the driving force and the evolution of the internal variables via
  \begin{equation}\label{eq:Dis_Ehat}
 	\mathcal{D}=\hat{\Ev}\cdot\dot{\Pv}_{I}\geq 0.
 \end{equation}
The constitutive minimization principle \eqref{eq:min_Pi_rev} may be extended to the case of irreversible processes. 
We introduce, according to \cite{MieheRosatoKiefer:2011} for the case of dissipative (or irreversible) processes, the potential $\Pi^{\tau}_{D}$ with its algorithmic expression 
\begin{align}
\begin{split}
	\Pi^{\tau}_{D}=&{\Psi}\left(\epst,\Dv,\Pv_{I}\right)-{\Psi}\left(\epst_n,\Dv_n,\Pv_{I,n}\right)\\
	&+D^{\tau}-\rev{W^\tau_{ext}}.
\end{split}
\end{align}
During the time interval $\tau$ the  energy stored is ${\Psi}\left(\epst,\Dv,\Pv_{I}\right)-{\Psi}\left(\epst_n,\Dv_n,\Pv_{I,n}\right)$  and $D^{\tau}$ is dissipated.
\par 
Analogously to (\ref{eq:min_Pi_rev}), for given external loads and states at time $t_n$, the strain, dielectric displacement and irreversible polarization satisfy the constitutive minimization principle which reads for dissipative processes
\begin{equation}\label{eq:min_Pi_ir}
	\begin{split}
	\left(\epst,\Dv,\Pv_{I} \right) =&\operatorname{Arg} \left\{\min_{\epst}\min_{\Dv}\min_{\Pv_{I}} ~\Pi^{\tau}_{D}\right\}
	\\ &~ s.t. \begin{cases}
		f_P(\hat{\Ev})\leq 0,\\
		f_S({\Pv_{I}})\leq 0,\\
		\opdiv\Dv=0.
	\end{cases}		
	\end{split}	
\end{equation} 
Note that, contrary to \eqref{eq:min_Pi_rev}, the minimization prinziple \eqref{eq:min_Pi_ir} for dissipative processes is constrained by inequality constraints. This has to be taken into account when solving for the free variables. \par 
\par 
As in the reversible case, an equivalent enthalpy-based formulation can be found. We define the total enthalpy $H$ as in \eqref{eq:H_int}, \eqref{eq:legendre_trans} but using the energy density $\psi$,
\begin{equation}\label{eq:h_min_dis}
	h\left(\sigmat,\Ev,\Pv_{I}\right)=\min_{\epst,\Dv}~\left(\psi\left(\epst,\Dv,\Pv_{I}\right)-\sigmat:\epst-\Dv\cdot\Ev\right).
\end{equation}
With the Lagrangian $\mathcal{L}$ and its potential $\rev{\Pi^{\tau}_{\mathcal{L}}}$ defined in analogy to \eqref{eq:H_ext}, \eqref{eq:Pi_tau_L} the following saddle point problem is obtained
\begin{equation}\label{eq:min_Entalpy_gen}
	\begin{split}
		\left(\sigmat,\Ev,\uv,\Pv_{I}\right)=&\operatorname{Arg} \left\{\max_{\Ev=-\nabla\phi}\max_{\sigmat}\min_{\uv}\min_{\Pv_{I}} ~\rev{\Pi^{\tau}_{\mathcal{L}}}\right\}\\
				\\ &~ s.t. \begin{cases}
			f_P(\hat{\Ev})\leq 0\rev{,}\\
			f_S(\rev{\Pv_{I}})\leq 0.\\
			\end{cases}	
	\end{split}
\end{equation}

\section{Incremental formulation of ferroelectricity} 
\label{sec:energy}
The focus of this section will be the realization of the energy density $\psi$ and the reformulation of the switching constraint $f_P({\hat{\Ev}})\leq  0$. The energy density will be chosen such that it reflects the constitutive equations \eqref{eq:const1}, \eqref{eq:const2} as well as the saturation condition \eqref{eq:sat}. 

We provide energy densities matching the constitutive equations from the previous section as a basis for all further deductions.
According to the literature (see e.g. \cite{Landis:2002,Kamlah:2001,Tichy:2010}), we follow the approach of additive decomposition of free energy. It decomposes into two parts,

\begin{eqnarray}\label{eq:free_e_H}
	\psi=\psi^R+ \psi^I,
\end{eqnarray}
with the reversible or stored part of the free energy $\psi^R$ and the additional contribution of the irreversible quantities to the free energy $\psi^I$. The reversible part $\psi^R$ is given by 
\begin{equation}\label{eq:PsiR}
	\begin{split}
	\psi^R=&\frac{1}{2} \Dv_R\cdot\left(\epsilont^{\epst}\right)^{-1}\cdot\Dv_R+\Dv_R\cdot\mathbf{h}:\epst_R\\
	&+\frac{1}{2}\epst_R:\Ct^{\Dv}:\epst_R,
	\end{split}
\end{equation}
\rev{with $\mathbf{h}$, the tensor of piezoelectric constants, $\epsilont^{\epst}$ the dielectric tensor at constant strain and $\Ct^{\Dv}$ the tensor of mechanical stiffness at constant dielectric displacement. These tensors are related to the material tensors in \eqref{eq:const1} in the standard way. 
\begin{align}
	&\St^{\Dv}(\Pv_{I})=\St^{\Et}-\dt(\Pv_{I})^T\cdot\betat^{\sigmat}\cdot\dt(\Pv_{I}),\\
	&\betat^{\epst}=\betat^{\sigmat}+\betat^{\sigmat}\cdot\dt(\Pv_{I}):\St^{\Dv}(\Pv_{I}):\dt(\Pv_{I})^T \cdot\betat^{\sigmat}),\\
	&\mathbf{h}(\Pv_{I})=\betat^{\epst}\cdot\dt(\Pv):\Ct^{\Ev},	\\
	&\Ct^{\Dv}(\Pv_{I})=\Ct^{\Ev}+\Ct^{\Ev}:\dt(\Pv_{I})\cdot\betat^{\epst}  \cdot\dt(\Pv_{I})^T   :\Ct^{\Ev},
\end{align}
with $\betat^{\sigmat}=(\epsilont^{\sigmat})^{-1}$ and $\betat^{\epst}=(\epsilont^{\epst})^{-1}$ and the tensor of mechanical compliance at constant dielectric displacement  $\St^{\Dv}$. This implies a non-trivial dependence on $\Pv_{I}$ for $\epsilont^{\epst}$, $\Ct^{\Dv}$ and $\mathbf{h}$. Note that, in the enthalpy-based formulations proposed in this work only much simpler structured material moduli $\epst^{\sigmat}$, $\St^{\Ev}$ and $\dt$ are needed. Of course non-trivial dependences on $\Pv_{I}$ may be avoided when assuming $\epsilont^{\epst}$ and $\Ct^{\Dv}$ to be constant and isotropic.   	}
\par
The additional contribution of the free energy is given by
\begin{equation}\label{eq:PsiI_el}
	\psi^I=\frac{1}{2}c\Pv_I\cdot\Pv_I,
\end{equation}
in the simplest case, where $c$ is a hardening parameter.
\par 
An additional term representing saturation is added in the following.
The saturation condition $f_S \leq 0$ (cmp. (\ref{eq:sat})) is reformulated as a variational inequality, involving the Lagrangian multiplier $\lambda_S$,
\begin{equation}\label{eq:sat_energy}
	\lambda_S f_S(\Pv_I)\leq 0 \qquad  \forall ~ \lambda_S \geq 0.
\end{equation}
To ensure the saturation condition, we add the following supremum term to the energy function:
\begin{equation}\label{eq:Psi_S}
	\psi=\psi^R+\psi^I+ \sup_{\lambda_S\geq 0} \lambda_S f_S(\Pv_I).
\end{equation}
For states of irreversible polarization $\Pv_{I}$ not violating the saturation condition ($f_S\leq 0$), the supremum in (\ref{eq:Psi_S}) takes zero value, and $\psi$ is not altered in comparison to (\ref{eq:free_e_H}). Otherwise, if $|\Pv_{I}|$ exceeds the saturation polarization, the  saturation condition is violated ($f_S>0$), and \rev{therefore}  $ \sup_{\lambda_S\geq 0} \lambda_S f_S$ tends to infinity, and with it the energy function $\psi$. Note, that for all admissible values of $\Pv_{I}$ - not violating the saturation condition - $\psi$ takes finite values, while for inadmissible {(improper)} values of $\Pv_{I}$,  the free energy function $\psi$ tends to infinity. In the spirit of optimization problems, this leads to solutions with admissible states of irreversible polarization $\Pv_{I}$. 
\par 
For further deductions we treat the Lagrange multiplier $\lambda_{S}$ as \rev{a} free variable. This is indicated by the index $S$, as we use
\begin{equation}\label{eq:psi_S}
	\psi^S=\psi^R+\psi^I+\lambda_{S}f_S
\end{equation}
In the final optimization problem $\lambda_{S}$ will be maximized.

\subsection{Transformation to enthalpy}
The enthalpy is a function of the free variables electric field $\Ev$,  mechanical stress $\sigmat$ and irreversible polarization $\Pv_I$ and is related to the free energy via Legendre transformation \eqref{eq:h_min_dis}. The enthalpy density $h^S$ including saturation condition \eqref{eq:psi_S} reads
\begin{equation}\label{eq:H}
\begin{split}
	h^S=&-\frac{1}{2}\Ev\cdot\epsilont^{\sigmat}\cdot\Ev-\frac{1}{2}\sigmat:\St^E:\sigmat-\sigmat:\dt\cdot\Ev\\
	&+\frac{1}{2}c\Pv_{I}\cdot\Pv_{I}-\Pv_{I}\cdot\Ev-\epst_I:\sigmat
	\\ &+\lambda_S f_S\left(\Pv_I\right).
	\end{split}
\end{equation} 
In the enthalpy setting the driving force is defined as the negative derivative of the total enthalpy $h^S$ by the irreversible quantities, here by the irreversible polarization $\Pv_I$. In contrast to the original model by Kamlah, the driving force $\hat{\Ev}$ contains saturation via the corresponding Lagrangian multiplier $\lambda_S$. The driving force reads
\begin{equation}\label{eq:defEhat}
\begin{split}
\hat{\Ev}=&-\frac{\partial h^S}{\partial\Pv_I} = \Ev + \sigmat:\left(\frac{\opd}{\opd \Pv_I}\dt\right)\cdot\Ev\\
&+ \sigmat:\left(\frac{\opd}{\opd\Pv_I}\epst_I\right) - c\Pv_I -\frac{\opd f_S\left(\Pv_{I}\right)}{\opd \Pv_I}\lambda_S.
\end{split}
\end{equation}
Note that \eqref{eq:defEhat} implies that via \eqref{eq:Dis_Ehat} the saturation is part of the dissipation $\mathcal{D}$.  
\subsection{Incremental optimization principle}

Again we consider the (finite) time interval \mbox{$\tau=t-t_n$}, with suffix $n$ denoting quantities at time $t_n$. For the modified enthalpy $h^S$ we proceed to a Lagrangian including the equilibrium condition in the same way as \eqref{eq:min_Pi_tau_H}
\begin{equation}\label{eq:ext_H_sat}
	\mathcal{L}^S=\int_\Omega h^S - \left(\opdiv\sigmat+\fv\right)\cdot\uv \operatorname{d\Omega}.
\end{equation}
The algorithmic expression of the corresponding potential $\Pi^{\tau}_S$ reads
\begin{align}\label{eq:pot_H}
\begin{split}
\rev{\Pi^{\tau}_S}=&\mathcal{L}^S(\sigmat,\Et,\uv,\Pv_{I},\lambda_{S})\\
&-\mathcal{L}^S(\sigmat_n,\Et_n,\uv_n,\Pv_{I,n},\lambda_{S,n})+D^{\tau}.
\end{split}
\end{align}
The saddle point problem reads
\begin{align}\label{eq:min_prinzip}
	\begin{split}
	 (\Ev,  &\sigmat, \uv, \Pv_{I},\lambda_{S} ) =\\
	   &\operatorname{Arg}\left\{\ \max_{\Ev=-\nabla\phi}~\max_{\sigmat}~\min_{\uv}~\min_{\Pv_I}~\max_{\lambda_S\geq 0} ~ \rev{\Pi^{\tau}_S} \right\}\\ 
	   &\text{s. t.} ~
		f_P\left(\hat{\Ev}\right)\leq 0. 
	\end{split}
\end{align}
In order to consider the switching condition, another Lagrangian multiplier $\lambda_{P}$ is introduced in analogy to \eqref{eq:sat_energy} via
 \begin{align}
 	\lambda_P f_P\left(\hat{\Ev}\right)&\leq 0  & \forall \lambda_P &\geq 0. \label{eq:KKT_pol}
 \end{align}
Consequently the Lagrangian $\mathcal{L}^S$ and the corresponding potential $\rev{\Pi^{\tau}_S}$ have to be updated by the terms in \eqref{eq:KKT_pol}. The Lagrangian $\mathcal{L}^P$, involving polarization, reads 
\begin{equation}\label{eq:H_ext_Pol}
	\mathcal{L}^P=\int_\Omega h^S + \left(\opdiv\sigmat+\fv\right)\cdot\uv  + \lambda_P\,f_P \operatorname{d\Omega}.
\end{equation}
The algebraic form of the corresponding potential $\rev{\Pi^{\tau}_P}$ is given by
\begin{align}\label{eq:pi_P}
\begin{split}
	\rev{\Pi^{\tau}_P}=&\mathcal{L}^P(\sigmat,\Et,\uv,\Pv_{I},\lambda_{S} \rev{,\lambda_{P}})-\\
	&\mathcal{L}^P(\sigmat_n,\Et_n,\uv_n,\Pv_{I,n},\lambda_{S} \rev{,\lambda_{P,n}})+ D_{\tau},
	\end{split}
\end{align} 
the corresponding saddle point problem, extended by the (independent) variable $\lambda_{P}$, reads 
\begin{equation}\label{eq:min_prinzip2}
\begin{split}
   &\left(\uv, \Ev, \sigmat,\Pv_{I},\lambda_{P},\lambda_{S} \right) =\\
& \operatorname{Arg}\left\{\min_{\uv}~\max_{\Ev=-\nabla\phi}~ \max_{\sigmat}~\min_{\Pv_I} \max_{\lambda_P\geq 0}~\max_{\lambda_S\geq 0}~\rev{\Pi^{\tau}_P} \right\}.
\end{split}
\end{equation}

\subsection{Variation of Optimization Problem}

{The saddle point problem \eqref{eq:min_prinzip2} may be solved by the method of variation. } As the Lagrangian multipliers $\lambda_S$ and $\lambda_P$ are constrained to be non-negative, the  variation of the potential $\rev{\Pi^{\tau}_P}$ leads to a variational inequality. Variation with respect to $\uv$, $\Ev$, $\sigmat$ and $\Pv_I$ can be done in the standard way. Arbitrary virtual values  $\delta\uv$, $\delta\sigmat$, $\delta\Ev$ and $\delta\Pv_I$ are admissible, only the essential (Dirichlet) boundary conditions on the corresponding boundaries $\Gamma_1$, $\Gamma_2$, $\Gamma_3$ and $\Gamma_4$ have to be fulfilled, respectively. A variational equation is obtained via
\begin{align}
  \frac{\partial\rev{\Pi^{\tau}_P}}{\partial\Ev}\cdot\delta\Ev
+ \frac{\partial\rev{\Pi^{\tau}_P}}{\partial\sigmat}:\delta\sigmat
+\frac{\partial\rev{\Pi^{\tau}_P}}{\partial\Pv_I}\cdot\delta\Pv_I
+\frac{\partial\rev{\Pi^{\tau}_P}}{\partial\uv}\cdot\delta\uv
&= 0.
\end{align}
  Note, that external (body) forces are already considered within the potential $\rev{\Pi^{\tau}_P}$.
 The saturation and polarization condition, which are given as inequalities, cannot be treated in the standard way. As the set of admissible $\lambda_P$ and $\lambda_S$ is restricted by inequalities $\delta\lambda_P$ and $\delta\lambda_S$ are not arbitrary. 
These conditions directly lead to variational inequalities, we refer to  the monograph \cite{HanReddy:1999} for an introduction into variational inequalities in the framework of elasto-plasticity. For both, $\lambda_P$ and $\lambda_S$ a test function $ \bar{\lambda}_P$ and $\bar{\lambda}_S$ is introduced. The potential $\rev{\Pi^{\tau}_P}$ is maximized with respect to $\lambda_P$  if and only if for all other admissible choices $\bar \lambda_P \geq 0$  there holds
\begin{align}
\bar \lambda_P\, f_P(\hat{\Ev}) &\leq \lambda_P\, f_P(\hat{\Ev}) & \forall \bar{\lambda}_P \geq 0 \label{eq:KKT_pol2}.
\end{align}
For the saturation and the corresponding Lagrangian multiplier $\lambda_S$, the situation is a bit more involved, as the driving force $\hat{\Ev}$ depends on $\lambda_S$. Writing 
\begin{equation}
	\bar{\hat \Ev} = \hat \Ev(\bar \lambda_S)
\end{equation}
for some  admissible test function $\bar \lambda_S \geq 0$, one finds that $\rev{\Pi^{\tau}_P}$ is maximized with respect to $\bar \lambda_S$ if and only if
\begin{align}
\begin{split}
\bar{\hat \Ev} \cdot \dot \Pv_I + \bar \lambda_S\, f_S\left(\Pv_I\right) + \lambda_P f_P(\bar{\hat \Ev}) \leq \\
 {\hat \Ev} \cdot \dot \Pv_I +  \lambda_S\, f_S\left(\Pv_I\right) + \lambda_P f_P({\hat \Ev})  \\ \forall \bar{\lambda}_S \geq 0 \label{eq:KKT_sat2}.
\end{split}
\end{align}
Inserting the definition of $\hat \Ev$ \eqref{eq:defEhat}, which implies that 
\begin{equation}
\bar{\hat \Ev} - \hat \Ev = -(\bar \lambda_S - \lambda_S) \frac{\opd f_S}{\opd\Pv_I},
\end{equation}
 the above inequality can  further be reduced to
\begin{align}
\begin{split}
(\bar \lambda_S - \lambda_S) \left( - \frac{\rev{\opd} f_S}{\rev{\opd} \Pv_I} \cdot \dot \Pv_I + f_S - \lambda_P \frac{{\opd}f_P}{{\opd}\hat{\Ev}} \frac{\opd f_S}{\opd \Pv_I} \right) \leq0 \\ \forall \bar{\lambda}_S \geq 0 \label{eq:KKT_sat22}.
\end{split}
\end{align}
Summing up, the variational inequality to hold for all admissible $\delta \Ev$, $\delta \sigmat$, $\delta \Pv_I$ and $\delta \uv$ as well as for all test functions $\bar \lambda_S\geq0$ and $\bar \lambda_P\geq0$ reads
\begin{align}\label{eq:var_ineq1}
\begin{split}
  \frac{\partial\rev{\Pi^{\tau}_P}}{\partial\Ev}\cdot\delta\Ev+
  + \frac{\partial\rev{\Pi^{\tau}_P}}{\partial\sigmat}:\delta\sigmat+\\
   \frac{\partial\rev{\Pi^{\tau}_P}}{\partial\Pv_I}\cdot\delta\Pv_I+  \frac{\partial\rev{\Pi^{\tau}_P}}{\partial\uv}\cdot\delta\uv+\\
    (\bar\lambda_S - \lambda_S)\,\left( f_S -  \frac{\opd f_S}{\opd \Pv_I} \cdot( \lambda_P\frac{\opd f_P}{\opd \hat \Ev} + \dot \Pv_I) \right)+&\\
     (\bar\lambda_P - \lambda_P)\,f_P      \leq &0.
    \end{split}
\end{align}

\subsection{Interpretation of the Lagrangian multipliers}
In this section we explain the meaning of the Lagrangian multipliers $\lambda_P$ and $\lambda_S$. For reasons of simplicity we restrict ourselves to the uncoupled, purely electrical problem. Note, that the conclusions of this section may be extended directly to the fully coupled problem, as the Lagrangian multipliers $\lambda_{P}$ and $\lambda_{S}$, as well as the corresponding inequalities read the same as for the mechanical problem with stress states $\sigmat=\nullt$. 

For this electric case, as all mechanical quantities vanish, the driving force reduces to the simple form
\begin{equation}
	\hat{\Ev}=\Ev-c\Pv_I-\lambda_S\frac{\Pv_I}{|\Pv_I|}.\label{eq:eHat_el}
\end{equation}
The variational inequality \eqref{eq:var_ineq1} reads
\begin{align}
\begin{split}
	\delta\left[-\frac{1}{2}\Ev\cdot\epsilont^{\rev{\sigmat}}\cdot\Ev - \Ev\cdot\Pv_I + \frac{c}{2}\Pv_I\cdot\Pv_I +\hat{\Ev}\cdot\dot{\Pv}_I\right]+& \\ 
		\lambda_P\left(\frac{\opd f_P}{\opd \hat{\Ev}}  \frac{\partial \hat{\Ev}}{\partial \Ev} \delta\Ev+\frac{\partial f_P}{\partial \Pv_I}\delta\Pv_I
		\right) +&\\
	 (\bar\lambda_S - \lambda_S)\left( f_S -  \frac{\opd f_S}{\opd \Pv_I} \cdot( \lambda_P\frac{\opd f_P}{\opd \hat \Ev} + \dot \Pv_I) \right) +&\\
	 (\bar\lambda_P - \lambda_P) f_P ~\leq &0. \label{eq:varineq_el2}
\end{split}
\end{align}

We consider an update step with the (updated) solutions $\Ev$, $\Pv_I$, $\lambda_{P}$ and $\lambda_{S}$. We introduce $\Delta\Pv_I$ as the update of the irreversible polarization between two calculation steps. As the considered process is rate independent and quasi-static the dissipation is given via $\hat{\Ev}\cdot\Delta\Pv_I$. Splitting up \eqref{eq:varineq_el2} into individual equations and inequalities for the single variational quantities we get

\begin{eqnarray}
	\delta\Ev\left(\rev{-}\epsilont^{\rev{\sigmat}}\cdot\Ev-\Pv_I+ \frac{\partial\hat{\Ev}}{\partial\Ev}\Delta\Pv_I + 
	\lambda_P\frac{\opd f_P}{\opd \hat{\Ev}}\frac{\partial\hat{\Ev}}{\partial\Ev} \right) &=0\\
	\delta\Pv_I \bigg( -\Ev +c\Pv_I  +\hat{\Ev} + \frac{\partial\hat{\Ev}}{\partial\Pv_I}\Delta\Pv_I+~& \nonumber\\
	\lambda_P\frac{\opd f_P}{\opd \hat{\Ev}}\frac{\partial\hat{\Ev}}{\partial\Pv_I} +  \lambda_S\frac{\opd f_S}{\opd\Pv_I} \bigg) &=0 \label{eq:varel_Pi2}\\
	(\bar\lambda_P - \lambda_P)\,f_P &\leq 0 \label{eq:varel_lambdaP}\\
  (\bar\lambda_S - \lambda_S)\,\left( f_S -  \frac{\opd f_S}{\opd \Pv_I} \cdot( \lambda_P\frac{\rev{\opd}f_P}{\rev{\opd} \hat \Ev} + \Delta \Pv_I) \right) &\leq 0. \label{eq:varel_lambdaS}
\end{eqnarray}

To get an interpretation of $\lambda_P$, consider a state, at which the polarization process has started, such that $f_P = 0$, $\hat{\Ev}\not\neq \nullv$ and $\lambda_P>0$ non zero. The irreversible polarization is further considered to be (and to stay) below the saturation level. In this case, the derivative of the switching condition is given by
 \begin{equation}\label{eq:dif_pol_cond_el}
	\frac{\rev{\opd} f_P}{\rev{\opd} \hat \Ev} = \hat \Ev / |\hat\Ev|.
\end{equation}
From the variation of the irreversible polarization $\delta\Pv_{I}$ in \eqref{eq:varel_Pi2}, taking into account the definition  (\ref{eq:eHat_el}) of $\hat{\Ev}$ and \eqref{eq:dif_pol_cond_el} we get a relation between $\lambda_P$ and $\Delta\Pv_I$ via
\begin{equation}\label{eq:interpret_lamP}
	\frac{\partial\hat{\Ev}}{\partial\Pv_I}\left(\Delta\Pv_I+\lambda_P\frac{\hat{\Ev}}{|\hat{\Ev}|}\right)=\nullv.
\end{equation}
From (\ref{eq:interpret_lamP}) one deduces that the polarization update $\Delta\Pv_I$ is in direction of $\hat{\Ev}$. As $\frac{\hat{\Ev}}{|\hat{\Ev}|}$ is a unit vector the norm of $\Delta\Pv_I$ equals $\lambda_P$. On the other hand, if the coercive field is not reached and $f_P < 0$, \eqref{eq:varel_lambdaP} ensures that $\lambda_P = 0$, and further $\Delta \Pv_I = \nullv$.

\begin{figure}[htp]
	\center
	\includegraphics[width=.35\textwidth]{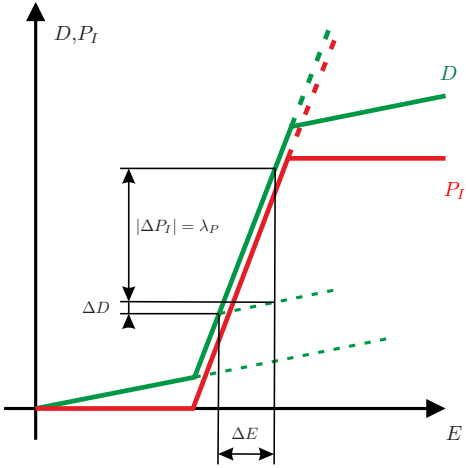}%
	\caption{Graphical interpretation of Lagrange multipliers for switching condition.} 
	\label{fig:lagpar_pol}
\end{figure}
\par 
Before we give an interpretation for $\lambda_S$, we first show that \eqref{eq:varel_lambdaS} really ensures the saturation condition. In case the coercive electric field is reached and $f_P = 0$, \eqref{eq:interpret_lamP} holds, and thereby the last term in \eqref{eq:varel_lambdaS} reduces to zero, leaving KKT conditions for the saturation condition only. On the other hand, if the coercive field is not reached and $f_P < 0$, we find $\lambda_P = 0$ and $\Delta \Pv_I = \nullv$, which again reduces the last term in \eqref{eq:varel_lambdaS} to zero. Thus the saturation condition has to be satisfied in any case.
\par 

\begin{figure}[htp]
	\includegraphics[width=.35\textwidth]{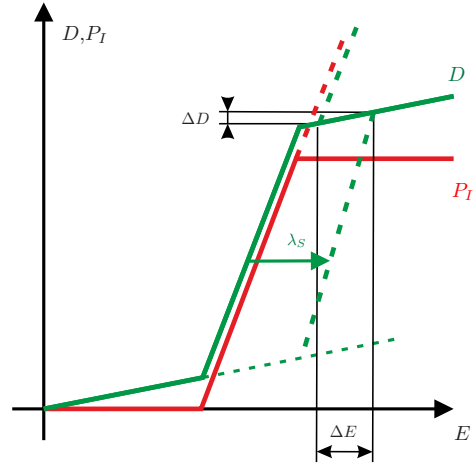}
	\caption{Graphical interpretation of Lagrange multipliers for  saturation condition.} 
	\label{fig:lagpar_sat}
\end{figure}

For the interpretation of $\lambda_S$ we consider a state where $f_S=0$, i.e. saturation polarization is reached. The driving electric field $\hat{\Ev}$ depending on $\lambda_{S}$ still has to satisfy the switching condition, 
\begin{align}
f_P = \left\|\Ev - c \Pv_I - \lambda_S \frac{\Pv_I}{|\Pv_I|}\right\| - E_C \leq 0.
\end{align}
The Lagrangian multiplier $\lambda_{S}$ grows with the electric field $\Ev$ as soon as further growth of the remanent polarization is prohibited by the saturation condition.
A graphical interpretation of $\lambda_P$ and  $\lambda_S$ can be found in Figure~\ref{fig:lagpar_pol} and Figure~\ref{fig:lagpar_sat}, respectively. 

\section{Finite element method }\label{sec:fem}

In this section, the proposed finite element discretization is described. 
For the electric potential continuous, nodal finite elements are used, while the polarization vector is approximated constant on each element. Also, the Lagrangian multipliers are realized taking one value per element.
As already done in the derivations of the previous sections, displacement and stress are considered independent unknowns. To do so, a mixed finite element scheme is used. The specific approach taken in this work uses tangential displacements and \rev{normal components of the normal stress vector $(\sigmat\cdot\nv)\cdot\nv$} as degrees of freedom, which motivate the abbreviation ``TDNNS''. This method was originally developed for elastic solids \rev{by} \cite{PechsteinSchoeberl:2011,PechsteinSchoeberl:2012} and later extended for linear piezoelectric materials \cite{PechsteinMeindlhumerHumer:2018,MeindlhumerPechstein:2018} and for geometrically nonlinear electro-mechanically coupled problems \cite{Pechstein:2019}. 
\rev{As the TDNNS method does not suffer from locking effects for elements of arbitrary aspect ratios, it is highly suitable for the discretization of thin structures.   }

This work is based on \cite{PechsteinMeindlhumerHumer:2018,MeindlhumerPechstein:2018} for linear \rev{piezoelectric materials}. 
Note that, for the mixed TDNNS method, the $\dt$-tensor formulation can be used directly, and there is no need to transfer the electric permittivity at constant stress $\epsilont^{\sigma}$ to that at constant strain $\epsilont^{\varepsilon}$ algebraically, as is necessary in standard formulations. As the underlying mathematics of the TDNNS method are rather involved, we skip it here, and only mention that work pairs such as $\int_\Omega \sigmat : \epst$ or $\int_\Omega\opdiv \sigmat \cdot \uv$ need to be considered in distributional sense and involve additional integration on element boundaries.

To solve the variational inequality, an active-set strategy is proposed. This means, in each iterative step, two active sets are identified beforehand: the switching active set $A_P$ where the switching condition shall be enforced as $f_P = 0$, and the saturation active set $A_S$ where the polarization condition is enforced as $f_S = 0$. On the other hand, for all elements not in the active set, the Lagrangian multipliers $\lambda_P$ and $\lambda_S$ are set to zero, respectively. For these fixed active sets, the nonlinear equations are solved by Newton's method, with $\lambda_S$ and $\lambda_P$ free in their respective active sets. On convergence, the Kuhn-Tucker compatibility conditions are checked. All elements where either saturation or switching condition are violated are added to the \rev{respective} active set. On the other hand, all elements where a Lagrangian multiplier is found negative are removed from the respective active set. An according algorithm can be found in Algorithm~\ref{alg:one}.
\par
\rev{Note that, in order to avoid numerical instabilities Kuhn-Tucker compatibility conditions in Algorithm~\ref{alg:one} are not checked for zero. The switching and saturation and saturation condition as well as the corresponding Lagrangian multipliers are compared to  a sufficiently  small value $\delta_P$ and $\delta_S$, respectively.
 These small offsets shall prevent recurring switching of elements between active and inactive state.
In our computations the values have been chosen \mbox{$\delta_P=E_{C}\cdot 10^{-3}$} and  \mbox{$\delta_S=P_{sat}\cdot 10^{-3}$}. These values lead to a maximum of three loops per load step for the three dimensional examples presented in the next chapter.    }

\subsection{Implementation of Dissipation}
In order to obtain an incremental formulation the dissipation has to be reformulated. The rate of irreversible polarization in the time interval $\tau=t-t_n$ is given by
\begin{equation}\label{eq:ir_pol_rate}
\dot{\Pv}_{I,n}=\left(\Pv_I-\Pv_{I,n}\right)/\tau=\Delta \Pv_{I}/\tau.
\end{equation}
Note that, in \eqref{eq:ir_pol_rate} a constant rate of polarization is assumed for one time interval. Taking into account \eqref{eq:Dis_Ehat} the time integral in \eqref{eq:dissipaton} is replaced by multiplication by $\tau$. The dissipated energy $D^{\tau} $ in the time interval $\tau$ is given by

\begin{equation}\label{eq:impl_ir_pol_rate}
	D^{\tau}=\int_\Omega \hat{\Ev}\cdot\left(\Pv_I-\Pv_{I,n}\right)\operatorname{d\Omega}=\int_\Omega \hat{\Ev}\cdot\Delta\Pv_I\operatorname{d\Omega}.
\end{equation}
Note, that in \eqref{eq:impl_ir_pol_rate} the driving force $\hat{\Ev}$ is taken at time $t$. This implies the backward Euler method for the driving force.  
\par 

\begin{algorithm}
	\KwData{values for $\uv_0, \sigmat_0, \phi_0, \Pv_{I0}, \lambda_{S0}, \lambda_{P0}$ and active sets $A_{S0}, A_{P0}$ from the last converged time step}
	\KwResult{values for  $\uv, \sigmat, \phi, \Pv_{I}, \lambda_{S}, \lambda_{P}$, active sets $A_S, A_P$}
	initialize unknows with initial data, active sets $A_S = A_{S0}, A_P = A_{P0}$\;
	\While{the active sets change}{
		solve non-linear problem with $\lambda_S = 0$ in $\Omega \backslash A_S$, $\lambda_P = 0$ in $\Omega \backslash A_P$ and $\lambda_S$ free in $A_S$, $\lambda_P$ free in $A_P$\;
		using starting values from last iteration\;
		\For{all elements $T$}{
			\uIf{ $T \in A_P, T \notin A_S$ and $\lambda_P < \rev{\delta_P}$}{ remove $T$ from $A_P$\;}
			\ElseIf{ $T \notin A_P$ and $f_P > \rev{\delta_P}$}{ add $T$ to $A_P$\;}
			\uIf{ $T \notin A_S, T \in A_P$ and $f_S > \rev{\delta_S}$}{ add $T$ to $A_S$\;}
			\ElseIf{ $T \in A_S$ and $\lambda_S < \rev{\delta_S}$}{ remove $T$ from $A_S$\;}
		}
	}
	\caption{Active set strategy for solving the variational inequality.}\label{alg:one}
\end{algorithm}

\section{NUMERICAL RESULTS}
In this section we show numerical results for several examples.
First we show the standard hysteresis curves for the fully coupled homogeneous (uni-axial) problem. We show that the formulation allows  mechanical depolarization. Our second example is the repolarization of ferroelastic material, that is initially polarized at a certain angle compared to the direction of applied electric field. Results are compared to the experiments carried out by \cite{HuberFleck:2001}. Next we show the mechanical depolarization of a fully polarized ferroelectric beam. This example was introduced as a ferroelastic benchmark in \cite{ZouariZinebBenjeddou:2011}. Finally we show the electric polarization of a bimorph structure. We use the same material data for all examples (except non-proportional loading) taken from \cite{ZouariZinebBenjeddou:2011} listed in Table \ref{tab:material}.

\begin{table}[htp]
	\caption{Material parameters\rev{.}}	
	\centering
	\begin{tabular}{ |l|l| } 
		\hline
		Young's modulus $Y$ & $10^4~\operatorname{MPa}$ \\ 
		Poisson's ratio $\nu$ & $0.3$ \\ 	
		Coercive field $E_C$ & $1~\operatorname{MV/m}$ \\ 
		Saturation strain $S_{sat}$& $0.002$\\ 
		Saturation polarization $P_{sat}$ & $0.3~\operatorname{C/m^2}$    \\  	
		Hardening parameter $c$ & $2\cdot10^6~\operatorname{Vm/C}$ \\	 
		Piezoelectric coefficient $d_p$ & $5.93\cdot10^{-10}~\operatorname{m/V}$\\ 
		Piezoelectric coefficient $d_n$ & $-2.74\cdot10^{-10}~\operatorname{m/V}$\\
		Piezoelectric coefficient $d_t$ & $7.41\cdot10^{-10}~\operatorname{m/V}$\\ 
		Permittivity $\epsilon$ & $1.5\cdot10^{-8}~\operatorname{C/Vm}$\\
		\hline
	\end{tabular}
	\label{tab:material}
\end{table}

\subsection{\rev{One}-dimensional loading}
\begin{figure}[htp]
	\centering
	\includegraphics[width=0.5\textwidth]{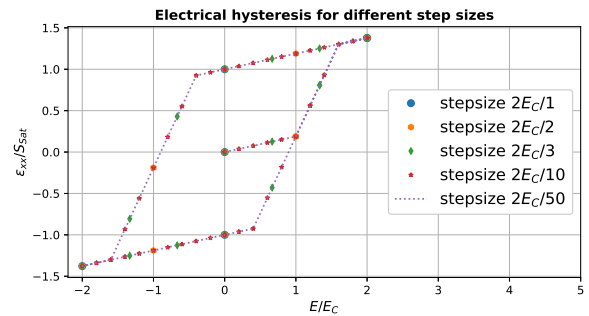}
	\caption{\rev{Electrical hysteresis of ferroelectric material.}}
	\label{fig:electric_polarization_el}
\end{figure}

\begin{figure}[htp]
	\centering
	\includegraphics[width=0.5\textwidth]{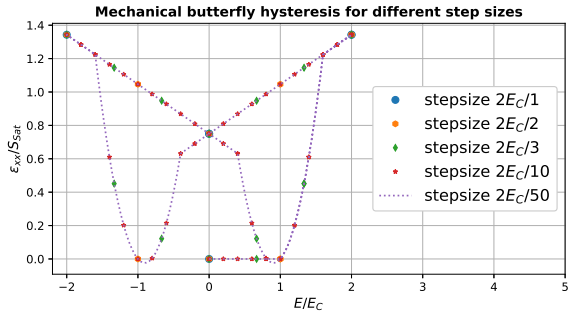}
	\caption{\rev{Mechanical butterfly hysteresis of ferroelectric material.}}
	\label{fig:electric_polarization_mech}	
\end{figure}

In our first example we  show that our implementation is capable of electrical polarization as well as electrical and mechanical depolarization. We consider an initially unpolarized unit square ($1~\operatorname{mm}$ by $1~\operatorname{mm}$) with electrodes at top and bottom. The material parameters are listed in Table \ref{tab:material}. Note, that for the two dimensional calculations plain strain is assumed. First the electrical polarization is shown. An electric field with a strength, for which the material will be fully polarized, here of two times the coercive field strength $2E_C$ is applied. Then the electric field is than reduced to $-2E_C$, and again increased up to $2E_C$. The resulting electric hysteresis, as well as the corresponding (mechanical) butterfly hysteresis are shown in Figure~\ref{fig:electric_polarization_el} and Figure~\ref{fig:electric_polarization_mech}, respectively. \rev{In order to show the convergence of the method, calculations are performed with different step sizes. Note that for the step size $2E_C$, full polarization and depolarization is reached within one load step, respectively. }   
\par 	

\begin{figure}
	\centering
	\includegraphics[width=0.5\textwidth]{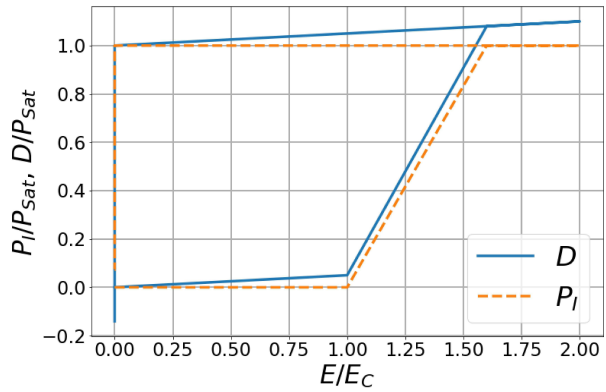}
	\caption{Electrical hysteresis - depolarization of ferroelectric material.}
	\label{fig:mech_polarization_el}
\end{figure}

\begin{figure}
	\centering
	\includegraphics[width=0.5\textwidth]{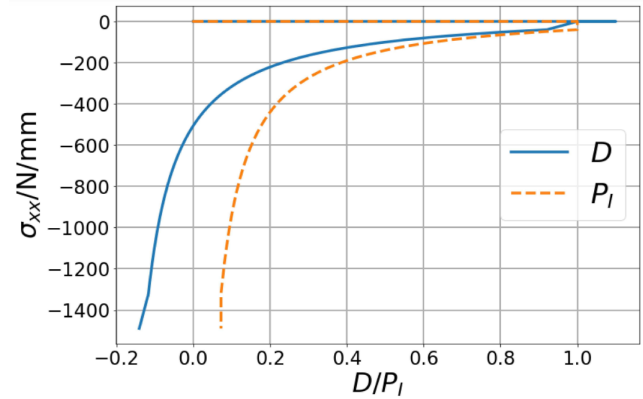}
	\caption{Mechanical butterfly hysteresis - depolarization of ferroelectric material.}
	\label{fig:mech_polarization_mech}
\end{figure}

Next we show the effect of mechanical depolarization. The specimen is electrically polarized (field strength $2 E_C$), then the electric field is reduced to $0$. A mechanical compressive stress (aligned in polarization direction) is applied to the polarized material and mechanical depolarization can be observed. The mechanical depolarization curve is shown in Figure~\ref{fig:mech_polarization_el} and Figure~\ref{fig:mech_polarization_mech}. 


\subsection{Non-proportional loading}

\begin{figure}[tbh]
	\centering
	\includegraphics[width=0.4\textwidth]{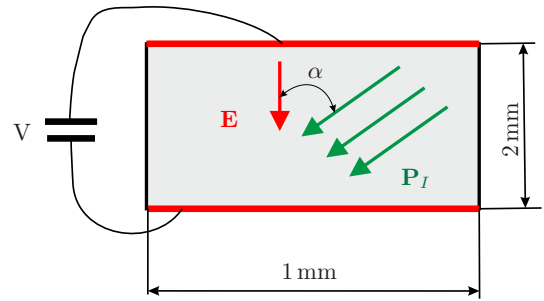}
	\caption{Sketch of repoling of ferroelectric material.}
	\label{fig:sketch_repol}
\end{figure}

In our second example we verify the capability of non-proportional loading. A rectangular specimen is initially fully poled in a certain direction. Electrodes are located at top and bottom, applying an electrical field in vertical direction. The direction of polarization, referring to the electric field, is denoted by the angle $\alpha$. The left and right surface of the specimen are \rev{insulated} and considered to be free of charges. All boundaries are free of mechanical stresses, only rigid body motion is restricted.  A sketch of the example, including geometry information, is shown in Figure~\ref{fig:sketch_repol}. The material data for for this example are taken form \cite{Kamlah:2001} (electrical parameters) and \cite{SemenovLiskowskyBalke:2010} (mechanical and coupling parameters). They are summarized  in Table~\ref{tab:material2}.
\par 
\begin{table}[htp]
	\caption{Material parameter\rev{s} for non-proportional loading. }	
	\centering
	\begin{tabular}{ |l|l| } 
		\hline
		Young's modulus $Y$ & $6.1\cdot10^4~\operatorname{MPa}$ \\ 
		Poisson's ratio $\nu$& $0.31$ \\ 	
		Coercive field $E_C$ & $0.8~\operatorname{MV/m}$ \\ 
		Saturation strain $S_{sat}$& $0.0046$\\ 
		Saturation polarization $P_{sat}$ & $0.23~\operatorname{C/m^2}$    \\  	
		Hardening parameter $c$ & $0.9\cdot10^6~\operatorname{Vm/C}$ \\	 
		Piezoelectric coefficient $d_p$ & $5.93\cdot10^{-10}~\operatorname{m/V}$\\ 
		Piezoelectric coefficient $d_n$ & $-2.74\cdot10^{-10}~\operatorname{m/V}$\\
		Piezoelectric coefficient $d_t$ & $7.41\cdot10^{-10}~\operatorname{m/V}$\\ 
		Permittivity $\epsilon$ & $6.2\cdot10^{-8}~\operatorname{C/Vm}$\\
		\hline
	\end{tabular}
	\label{tab:material2}
\end{table}

{As a consequence of the boundary conditions, which imply $\Dv\cdot\nv=0$ at \rev{insulated} boundaries, for arbitrary \rev{angle} $\alpha$, homogeneous polarization is an incompatible state. Here this issue is solved by  \rev{increasing} the polarization in small steps and then solving for zero voltage (both electrodes grounded). Values of dielectric displacement are taken in the  middle of the specimen, where boundary conditions do not affect the results.  A more involved discussion about the effect of boundary conditions can be found in \cite{StarkNeumeisterBalke:2016a}.}

\begin{figure}[thb]
	\centering
	\includegraphics[width=0.5\textwidth]{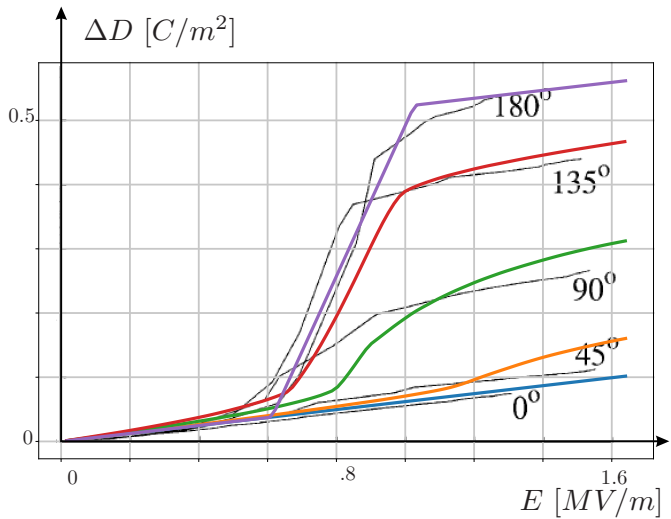}
	\caption{Non-proportional loading - repoling of ferroelectric material.}
	\label{fig:repol}
\end{figure}

  An electric field of size $2E_C$ is applied though the electrodes. Due to the applied electric field the direction of polarization is changed. The electric field and the change of dielectric displacement are measured at the center of the specimen. As a result the change of dielectric displacement in vertical direction over the applied electric field is shown for five particular \rev{angles} alpha, ($0^\circ$, $45^\circ$, $90^\circ$, $135^\circ$, $180^\circ$) in Figure~\ref{fig:repol}. The results are visually compared to those from the experiments in \cite{HuberFleck:2001}. The experimental data is shown in black lines in the back in Figure~\ref{fig:repol}, our results are colored. We find good correlation between the measurements and the calculations.
  
\subsection{Depolarization of a ferroelectric beam}
\begin{figure}[tph]
	\centering
	\includegraphics[width=0.4\textwidth]{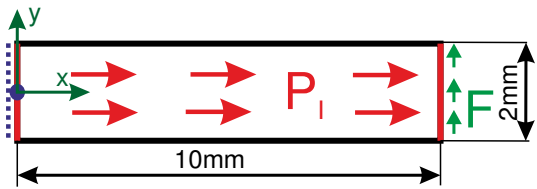}
	\caption{Sketch of fully polarized cantilever beam with tip load.}
	\label{fig:beamsketch}
\end{figure}

 Our third example is a benchmark for ferroelastic media taken from  \cite{ZouariZinebBenjeddou:2011}. It is a cantilevered beam, fully polarized in its longitudinal direction and loaded with a tip force at its end. A sketch of the problem can be found in Figure \ref{fig:beamsketch}. The computation is performed for a 3D-setting with square cross section ($2~\operatorname{mm}$ depth of the specimen). Grounded electrodes ($\phi=0$) are located at the clamped and at the tip face (illustrated in red). The clamping is realized by restricting longitudinal displacement \rev{and rigid body motions} at the clamped face (blue dashed line).
 \par
\begin{figure}[tph]
 	\centering
 	\includegraphics[width=0.4\textwidth]{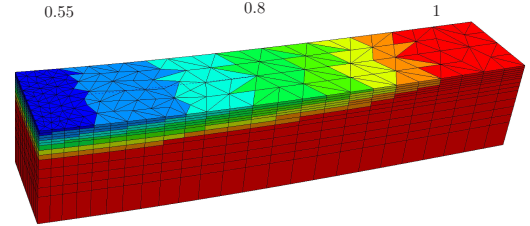}
 	\caption{Depolarization of the ferroelectric beam under tip load.}
 	\label{fig:beam_Pi}
 \end{figure}

\begin{figure}[tph]
	\centering
	\includegraphics[width=0.4\textwidth]{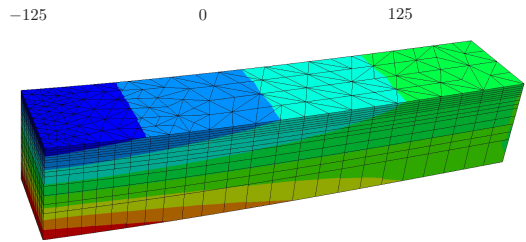}
	\caption{Stress $\sigma_{xx}$ of the ferroelectric beam under tip load.}
	\label{fig:beam_sxx}
 \end{figure}

 
The tip force is realized as stress boundary condition, applying $\sigma_{xy}=2 \operatorname{N/mm^2}$ at the tip face. Taking into account the dimensions, this load equals the tip force of $8\operatorname{N}$ in \cite{ZouariZinebBenjeddou:2011}. The applied load causes bending of the beam, which leads to a compressive stress at the top face and therefore mechanical depolarization can be observed. The resulting polarization is shown in Figure~\ref{fig:beam_Pi}. 
The resulting stresses are shown in Figure \ref{fig:beam_sxx}. Note that the irreversible polarization takes one (constant) value in each element, while the stress is interpolated.  Due to the unsymmetrical remanent straining an unsymmetrical distribution of the bending stresses is to be observed. 
\rev{In Figure~\ref{fig:beam_Pi_Sigma} the bending stress $\sigma_{xx}$ on top and bottom of the beam, as well as the irreversible polarization on the top are plotted over the length of the beam. 
} 
Note that, in contrast to the reference, only ferroelectric but no ferroelastic effects are considered in the constitutive model. 
\rev{While in \cite{ZouariZinebBenjeddou:2011} a decrease of the irreversible polarization of $0.82 {P_{sat}}$ is reported at the clamped face, here the minimum irreversible polarization  is $0.5485 {P_{sat}}$. Furthermore, in the reference the the top surface of the beam is fully polarized at locations $x>\SI{6}{\milli\meter}$, here the state of full polarization is present at positions $x>\SI{8.78}{\milli\meter}$. }

\rev{ \begin{figure}[tph]
 	\centering
 	\includegraphics[width=0.5\textwidth]{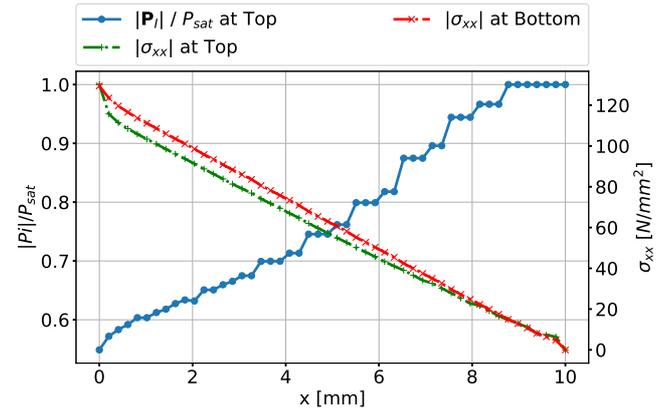}
 	\caption{\rev{Stress at top and bottom, polarization at top for polarized cantilever beam with tip load.}}
 	\label{fig:beam_Pi_Sigma}
 \end{figure}}

\subsection{Bimorph structure}

\begin{figure}[tph]
	\centering
	\includegraphics[width=0.4\textwidth]{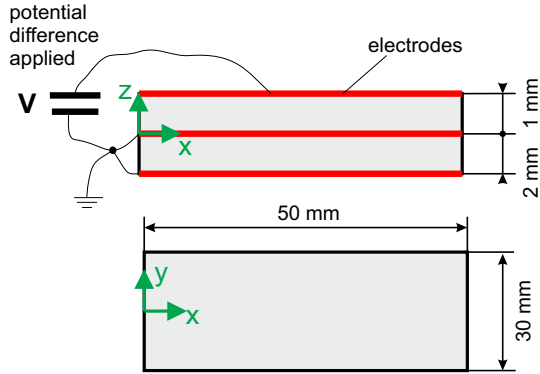}
	\caption{Sketch of bimorph structure.}
	\label{fig:sketch-bimorth}
\end{figure}

\begin{figure}[tph]
	\centering
	\includegraphics[width=0.4\textwidth]{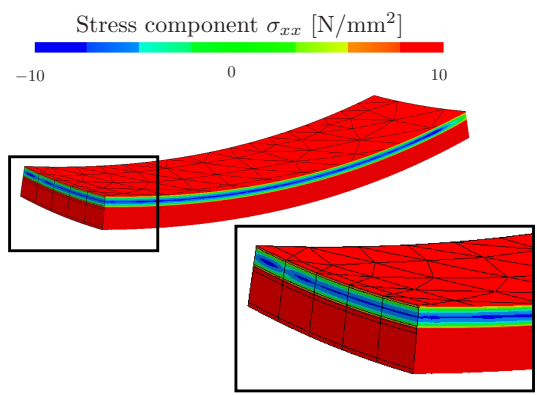}
	\caption{Stress component $\sigma_{xx}$ after polarization bimorph structure.}
	\label{fig:bimorth-S11}
\end{figure}

\begin{figure}[tph]
	\centering
	
	\includegraphics[width=0.4\textwidth]{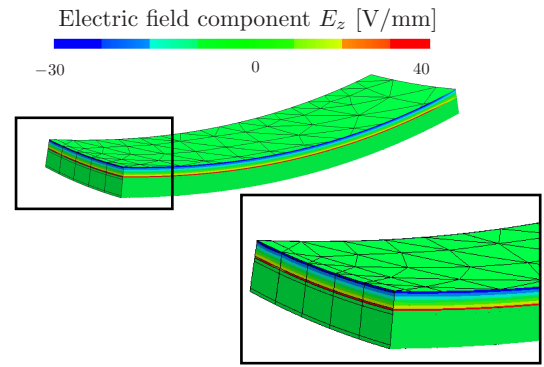}
	\caption{Electric field component $E_z$ after polarization bimorph structure.}
	\label{fig:bimorth-E3}
\end{figure}

\begin{figure}
	\centering
	\includegraphics[width=0.5\textwidth]{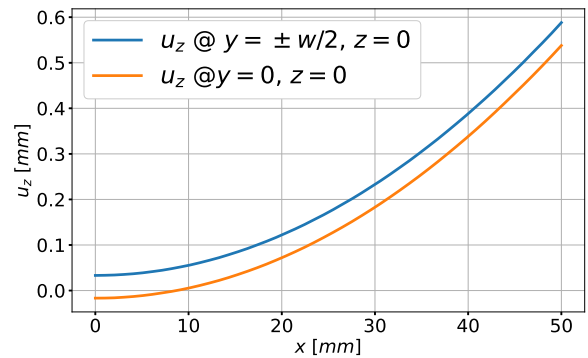}
	\caption[Deflection of the  bimorph structure]{\rev{Deflection of the  bimorph structure after polarization at the center line ($y=0$, $z=0$) and on the lateral surfaces ($y=\pm w/2$, $z=0$).}}
	\label{fig:bimorth-defl}
\end{figure}

In our last example we show the polarization of a bimorph structure. A sketch of the structure, including geometry parameters is shown in Figure~\ref{fig:sketch-bimorth}. The bimorph consists of two layers of different thickness, both of the same ferroelectric material with material parameters according to Table~\ref{tab:material}. Only the upper layer is electrically active, as the lowest and the middle electrode are grounded. The specimen is cantilevered in the same way as in the previous example, such that elongation within the fixed plane is not restricted. Figure~\ref{fig:sketch-bimorth} shows the  mesh used for the calculation. As the used TDNNS elements do not suffer \rev{from} locking effects, a very \rev{anisotropic} (spatial) discretization can be chosen. In thickness direction in total six mesh layers are used, with very thin layers near the electrodes. 
The upper layer is fully polarized by  applying a sufficient high electric field. After polarization the poling field is removed, all electrodes are grounded. Due to remanent straining only in the upper layer, the structure is deformed. Note that the final configuration is not free of mechanical stress, which further leads to (rather high) electric fields.
The resulting stress component $\sigma_{xx}$, as well as the  electric field in thickness direction $E_z$ are shown in Figure~\ref{fig:bimorth-S11} and Figure~\ref{fig:bimorth-E3}, respectively. The structure is shown in a deformed configuration with displacement scaled by a factor of $20$.
\rev{
	In Figure~\ref{fig:bimorth-defl} the deflection $u_z$ is evaluated along two horizontal lines, one of them the bimorph's center line ($y=0$, $z=0$), one on the lateral surface ($z=0$, $y=\pm w/2$). The deflection at the end of the center line is \mbox{$u_{z}(x=b,y=0,z=0)=\SI{0.5379}{\milli\meter}$} , the deflection of the corner points is \mbox{$u_{z}(x=b,y=\pm w/2,z=0)=\SI{0.5882}{\milli\meter}$}.
	
}
 
\begin{figure}[tph]
	\centering
	\includegraphics[width=0.4\textwidth]{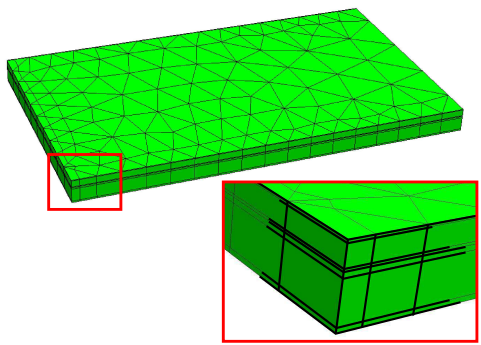}
	\caption{Mesh for calculation of bimorph structure.}
	\label{fig:mesh-bimorth}
\end{figure}

\section{Conclusion}
In this paper we have shown the implementation of (known) material laws in the framework of variational inequalities. This allows the direct implementation of inequalities in frameworks based on variational concepts as e.g. finite elements. As shown, inequality constraints can be included in the calculation process and are taken into account directly (no return mapping is used). A finite element discretization is chosen  using normal-normal stress and tangential displacements for the mechanical degrees of freedom, and electric potential and  remanent polarization for  the electrical degrees of freedom as well as two (additional) Lagrangian multipliers for the inequality constraints. The choice of free variables allows for direct use of the piezoelectric tensor $\dt$.    \par 
The  performance of the method was shown within several numerical examples, where results according to the literature were provided.

\begin{acks}
Martin Meindlhumer acknowledges support of Johannes Kepler University Linz, Linz Institute of Technology (LIT).
\end{acks}

\bibliographystyle{SageH}      
\bibliography{Polarization}


%
%
%

\end{document}